\documentclass[11pt, reqno]{amsart}
\usepackage{amsmath, amsthm, amscd, amsfonts, amssymb, graphicx, color}
\usepackage{bbm,mathtools}
\usepackage{booktabs}
\usepackage[bookmarksnumbered, colorlinks, plainpages]{hyperref}
\usepackage{enumitem}
\usepackage{mathrsfs}
\hypersetup{colorlinks=true,linkcolor=red, anchorcolor=green, citecolor=cyan, urlcolor=red, filecolor=magenta, pdftoolbar=true}

\newtheorem{theorem}{Theorem}[section]
\newtheorem{lemma}[theorem]{Lemma}
\newtheorem{proposition}[theorem]{Proposition}
\newtheorem{corollary}[theorem]{Corollary}
\theoremstyle{definition}
\newtheorem{definition}[theorem]{Definition}

\theoremstyle{remark}
\newtheorem{remark}[theorem]{Remark}
\newtheorem{claim}[theorem]{Claim}

\newcommand{\uinorm}[1]{{\left\vert\kern-0.25ex\left\vert\kern-0.25ex\left\vert #1
 \right\vert\kern-0.25ex\right\vert\kern-0.25ex\right\vert}}

\newcommand{\id}{\mathrm{id}}

\newcommand{\CP}{\mathrm{CP}}

\newcommand{\diag}{\mathrm{diag}}

\newcommand{\cconv}{\mathrm{cconv}}

\newcommand{\1}{\mathbbm{1}}

\begin{document}
\title[Matrix-Test Duality for $C^*$-Convex    Families of CP Maps]{Matrix-Test Duality: A Support-Function Characterization for $C^*$-Convex Families of CP Maps}

\author{Mohsen Kian$^1$}
\address{$^1$Department of Mathematics, University of Bojnord, P. O. Box 1339, Bojnord 94531, Iran}
\email{kian@ub.ac.ir }

\author{Mario Krni\' c$^2$}
\address{$^2$University of Zagreb, Faculty of Electrical Engineering and Computing, Unska
	3, 10000 Zagreb, Croatia}
\email{mario.krnic@fer.hr}


\begin{abstract}
We develop a matrix-test framework for $C^*$-convex families of
completely positive maps $\CP(\mathscr S,\mathscr T)$, where $\mathscr S$ is an
operator system and $\mathscr T$ is a unital $C^*$-algebra. Matrix tests
$(k,f,s)$ induce evaluation functionals $\Phi\mapsto f(\Phi_k(s))$ and generate
a natural weak topology $\tau$ on
$\mathcal E=\operatorname{span}_{\mathbb C}(\CP(\mathscr S,\mathscr T))$.

Our main result provides a support-function/separation characterization of the
$\tau$-closure $\overline{\cconv(\mathcal K)}^{\tau}$ of the $C^*$-convex
hull of a family $\mathcal K\subseteq\CP(\mathscr S,\mathscr T)$ in terms of
matrix-test inequalities. A key technical tool is a folding lemma that
represents every nonzero finite complex linear combination of test
functionals, up to a positive scalar factor, by a single matrix test at a
suitable block level. As consequences, we obtain a single-test witness for
non-membership, support-function criteria for inclusion and equality of the
corresponding $\tau$-closed hulls, and, under
$0\in\overline{\cconv(\mathcal K)}^{\tau}$, an exact normalized
bipolar-type reconstruction formula. We also show that $\tau$ is already
generated by level-$1$ tests.
\end{abstract}

\maketitle

\section{Introduction}

Duality and separation methods are fundamental tools in the study of positivity
and convexity in operator-algebraic settings. In classical  convex
analysis, closed convex sets are described by the continuous affine inequalities
they satisfy. In noncommutative settings, this principle survives, but the
relevant geometry is richer: one works with matrix levels, positivity
constraints, and structured notions of convexity. This perspective is central
in matrix-convex and tracial duality theories, including the matrix
Hahn--Banach/bipolar theory of Effros and Winkler
\cite{EffrosWinkler1997} and the tracial framework developed by Helton, Klep,
and McCullough \cite{HeltonKlepMcCullough2017}.

Matrix and noncommutative convexity already admit intrinsic
map-theoretic formulations through the matrix and nc state spaces of
operator systems.  In particular, Webster and Winkler showed that
compact matrix convex sets in locally convex spaces are, up to matrix
affine homeomorphism, matrix state spaces of operator systems
\cite[Proposition~3.5]{WebsterWinkler1999}; the noncommutative convexity framework of Davidson and Kennedy also includes infinite-dimensional levels \cite{DavidsonKennedy2025}.  Thus, the distinction relevant to the
present paper is not one between a point-based theory and a map-based
theory.

Our focus is instead on $C^*$-convex subsets of the fixed space
$\CP(\mathscr S,\mathscr T)$, where $\mathscr S$ is an operator system
and $\mathscr T$ is a prescribed unital $C^*$-algebra.  When
$\mathscr T=M_n$, the defining $C^*$-convex combinations are the
same-level matrix convex combinations; in the present setting, however,
no graded family or compatibility between different matrix levels is
part of the data.  We develop a support-function and separation
framework adapted to this fixed-codomain setting.
 
A complementary line of work concerns $C^*$-convexity and its geometry,
including structural and extremal aspects of $C^*$-convex sets
\cite{Farenick1993,kian2017, kian2026,LoeblPaulsen1981}.
Connections between $C^*$-convexity and completely positive maps have also been
studied, for example in \cite{Magajna2016}, and recent work has examined
$C^*$-extreme structure for classes of CP maps
\cite{Bala2024,BhatKumar2022}.

Against this background, we introduce a support-function framework
based on \emph{matrix tests}, adapted to the fixed-codomain setting.
A matrix test is a triple $(k,f,s)$, where $k\in\mathbb N$,
$f\in S(M_k(\mathscr T))$ is a state on $M_k(\mathscr T)$, and
$s\in M_k(\mathscr S)$. It evaluates a map
$\Phi\in\CP(\mathscr S,\mathscr T)$ by
\[
\langle (k,f,s),\Phi\rangle = f(\Phi_k(s)).
\]
Here the level $k$ belongs to the test rather than to a graded family
of maps.

Passing to the ambient vector space
\[
\mathcal E
=
\operatorname{span}_{\mathbb C}
\bigl(\CP(\mathscr S,\mathscr T)\bigr),
\]
the same formula defines a complex-linear functional on $\mathcal E$
for each matrix test. Letting $\mathcal F$ denote the complex-linear
span of these functionals, we obtain the weak topology
\[
\tau=\sigma(\mathcal E,\mathcal F).
\]

We give a support-function/separation characterization of the
$\tau$-closure
\[
\overline{\cconv(\mathcal K)}^{\,\tau}
\]
of the $C^*$-convex hull of a family
$\mathcal K\subseteq\CP(\mathscr S,\mathscr T)$. In concrete terms, it
says that a map
$\Phi_0\in\CP(\mathscr S,\mathscr T)$ belongs to
$\overline{\cconv(\mathcal K)}^{\,\tau}$ exactly when it satisfies the
matrix-test inequalities determined by the support function of
$\mathcal K$. The separation step itself is a standard Hahn--Banach
argument in the dual pair $(\mathcal E,\mathcal F)$. The matrix-test
framework and the folding mechanism represent every finite complex
linear combination of test functionals with coefficients not all zero
as a positive scalar multiple of a single test functional at a
suitable block level, thereby turning the abstract separation
criterion into concrete matrix-test inequalities.

The paper is organized as follows. In Section~2, we introduce matrix
tests, the associated weak topology $\tau$, and the polar operations
used later in the paper. We also show that $\tau$ is already generated
by level-$1$ tests. Section~3 contains the folding lemma and the
separation result. In Section~4, we prove the main support-function
characterization and derive the corresponding non-membership,
inclusion, equality, and bipolar-type criteria.


\section{Preliminaries and Setup}
This section fixes notation and collects the basic definitions used throughout
the paper.
In what follows, $\mathscr{S}$ is an operator system and
$\mathscr{T}$ is a unital $C^*$-algebra. We denote by $\mathscr{S}^+$ the cone of positive elements of $\mathscr{S}$. A linear mapping $\Phi:\mathscr{S}\to\mathscr{T}$ is called positive if $\Phi(\mathscr{S}^+)\subseteq \mathscr{T}^+$.

For $k\in\mathbb N$, the \emph{$k$-ampliation} of $\Phi$ is
\[
\Phi_k := \id_{M_k}\otimes \Phi :
M_k(\mathscr{S})\longrightarrow M_k(\mathscr{T}),\qquad
\Phi_k([s_{ij}])=[\Phi(s_{ij})].
\]
 If $\Phi_k $ is positive, then $\Phi$ is called $k$-positive. If $\Phi$ is $k$-positive for all $k\in\mathbb{N}$, then it is called completely positive. We denote by $\CP(\mathscr{S},\mathscr{T})$ the set of all completely positive
maps $\Phi:\mathscr{S}\to \mathscr{T}$.

We now recall the $C^*$-convex structure on $\CP(\mathscr{S},\mathscr{T})$ that
will be used throughout the paper.
Given $\Phi_1,\dots,\Phi_n\in \CP(\mathscr{S},\mathscr{T})$ and coefficients
$a_1,\dots,a_n\in \mathscr{T}$ satisfying
$\sum_{i=1}^n a_i^*a_i=\1_\mathscr{T}$, define
\begin{equation}\label{comb}
   \Phi(\,\cdot\,)
   :=\sum_{i=1}^n a_i^*\,\Phi_i(\,\cdot\,)\,a_i .
\end{equation}
Since each map $x\mapsto a_i^*xa_i$ on $\mathscr{T}$ is completely positive,
the map $\Phi$ in \eqref{comb} again belongs to $\CP(\mathscr{S},\mathscr{T})$.
We call \eqref{comb} a \emph{$C^*$-convex combination} of the maps $\Phi_i$.

A subset $\mathcal K\subseteq\CP(\mathscr{S},\mathscr{T})$ is said to be
\emph{$C^*$-convex} if it is closed under all finite combinations of the form
\eqref{comb}. We write $\cconv(\mathcal K)$ for the smallest $C^*$-convex
subset of $\CP(\mathscr{S},\mathscr{T})$ containing $\mathcal K$.

For each $k\ge1$, let
\[
S(M_k(\mathscr{T})) :=
   \{\, f:M_k(\mathscr{T})\to\mathbb C
   \mid f \text{ linear, positive, and } f(\1_{M_k(\mathscr{T})})=1 \,\}
\]
denote the state space of the unital $C^*$-algebra $M_k(\mathscr{T})$.
We refer to its elements as \emph{matrix states}. We also set
\[
\mathcal E:=\mathrm{span}_{\mathbb C}\bigl(\CP(\mathscr{S},\mathscr{T})\bigr),
\]
i.e.,  $\mathcal E$ is a complex vector space of linear maps
$\mathscr{S}\to\mathscr{T}$.

\begin{definition}
We define a \emph{matrix test} as  a triple $(k,f,s)$ consisting of
\begin{enumerate}
  \item a level $k\in\mathbb N$,
  \item a matrix state $f\in S(M_k(\mathscr{T}))$,
  \item a test element $s\in M_k(\mathscr{S})$.
\end{enumerate}

Let
\[
\mathfrak T:=\{(k,f,s)\mid k\in\mathbb N,\ f\in S(M_k(\mathscr T)),\
s\in M_k(\mathscr S)\}
\]
denote the collection of all matrix tests.
This is not a vector space (in particular, the matrix level $k$ varies).
Therefore, whenever we discuss convexity or closedness on the \emph{test side},
we do so after passing to the associated functionals
$\ell_{k,f,s}\in\mathcal F$.

For each matrix test $(k,f,s)$, define
\[
\ell_{k,f,s}:\mathcal E\to\mathbb C,\qquad
\ell_{k,f,s}(\Phi):=f\!\big(\Phi_k(s)\big).
\]
This is well defined because $\Phi\in\mathcal E$ is a linear map
$\mathscr{S}\to\mathscr{T}$, hence $\Phi_k(s)\in M_k(\mathscr{T})$, and
$f$ is a linear functional on $M_k(\mathscr{T})$.

For $\Phi\in\CP(\mathscr{S},\mathscr{T})$, we write the corresponding
\emph{evaluation pairing} as
\begin{equation}\label{eq:pairing}
   \big\langle (k,f,s),\,\Phi \big\rangle
      := \ell_{k,f,s}(\Phi)= f\!\big( \Phi_k(s) \big)\in\mathbb C .
\end{equation}
\end{definition}
The pairing~\eqref{eq:pairing} has the following operational interpretation:
one inputs $s\in M_k(\mathscr{S})$ into the amplified map $\Phi_k$, obtains
$\Phi_k(s)\in M_k(\mathscr{T})$, and then evaluates the output using the
matrix state $f$. The polar and support-function constructions collect the
matrix tests that impose uniform one-sided bounds on these evaluations.

\begin{lemma}\label{lem:pairing-basic}
Let $(k,f,s)$ be a matrix test, and let $\Phi\in\mathcal E$.
Then:
\begin{enumerate}
\item[(a)] The map
\[
   \Phi\longmapsto \ell_{k,f,s}(\Phi)=\langle (k,f,s),\Phi\rangle
\]
is a complex-linear functional on $\mathcal E$.

\item[(b)] If $\Phi\in\CP(\mathscr S,\mathscr T)$ and $s\ge 0$ in
$M_k(\mathscr S)$, then
\[
\langle (k,f,s),\Phi\rangle = f(\Phi_k(s))\ge 0.
\]

\item[(c)] For fixed $k$ and fixed $\Phi\in\mathcal E$, the map
$(f,s)\mapsto f(\Phi_k(s))$ is affine in $f\in S(M_k(\mathscr T))$
and complex-linear in $s\in M_k(\mathscr S)$.
\end{enumerate}
\end{lemma}

\begin{definition}\label{def:polars}
For a subset $\mathcal L\subseteq\CP(\mathscr{S},\mathscr{T})$, define
\begin{align}\label{eq:pol_K}
   \mathcal L^{\circ_{C^*}}
   :=\Bigl\{(k,f,s)\,\Big|\,
     \forall\,\Phi\in\mathcal L:
     \Re\, f\!\big(\Phi_k(s)\big)\le 1\Bigr\},
\end{align}
where $\Re z$ means the real part of the complex number $z$.

For a subset $\mathcal K\subseteq\CP(\mathscr{S},\mathscr{T})$, define its
\emph{$C^*$-convex polar} by
\begin{equation}\label{eq:cconv-polar}
   \mathcal K^{\circ}
   := (\cconv(\mathcal K))^{\circ_{C^*}}.
\end{equation}
Conversely, for a collection of matrix tests $\mathcal C$ we define
\begin{equation}\label{eq:polar-of-C}
   \mathcal C^{\circ_{C^*}}
   :=\Bigl\{\Phi\in\CP(\mathscr{S},\mathscr{T})\,\Big|\,
     \forall\,(k,f,s)\in\mathcal C:
       \Re\, f\!\big(\Phi_k(s)\big)\le 1\Bigr\}.
\end{equation}
\end{definition}

\begin{remark}\label{rem:one-sided-polar}
The polar convention in this paper is one-sided: it is defined by inequalities
of the form $\Re\,\ell_{k,f,s}(\Phi)\le 1$. This is the normalization naturally
adapted to convex separation (Hahn--Banach) in the dual pair
$(\mathcal E,\mathcal F)$ introduced below.

If $\mathcal K$ is already $C^*$-convex, then
\[
\mathcal K^\circ=\mathcal K^{\circ_{C^*}}.
\]
\end{remark}

Both polar constructions are monotone and order-reversing:
\[
\mathcal K_1\subseteq\mathcal K_2 \;\Rightarrow\;
   \mathcal K_2^{\circ}\subseteq\mathcal K_1^{\circ},
\qquad
\mathcal C_1\subseteq\mathcal C_2 \;\Rightarrow\;
   \mathcal C_2^{\circ_{C^*}}\subseteq\mathcal C_1^{\circ_{C^*}}.
\]

\subsection{The natural locally convex topology}

The evaluation functionals introduced above also determine the weak topology
used throughout the paper. In this subsection we study the corresponding
locally convex topology on $\mathcal E$ and its description in terms of the
span of matrix-test functionals. In particular, this topology is independent of
the one-sided normalization in Definition~\ref{def:polars}.

Define a family of seminorms on $\mathcal E$ by
\begin{equation}\label{eq:seminorm}
   p_{k,f,s}(\Phi) := \big| f\!\big(\Phi_k(s)\big)\big|
   = |\ell_{k,f,s}(\Phi)|.
\end{equation}
Let $\tau$ be the locally convex topology on $\mathcal E$ generated by all such
seminorms. Equivalently,
\[
\tau=\sigma(\mathcal E,\mathcal F),
\qquad
\mathcal F:=\mathrm{span}\{\ell_{k,f,s}\mid
k\in\mathbb N,\ f\in S(M_k(\mathscr{T})),\ s\in M_k(\mathscr{S})\}
\subseteq \mathcal E^* .
\]

\begin{lemma}\label{lem:tau-dual}
For each matrix test $(k,f,s)$, the functional $\ell_{k,f,s}$ is
$\tau$-continuous on $\mathcal E$. Conversely, every $\tau$-continuous
linear functional on $\mathcal E$ belongs to $\mathcal F$.
\end{lemma}

\begin{proof}
By construction, \(\tau=\sigma(\mathcal E,\mathcal F)\) is the weak (initial)
locally convex topology on \(\mathcal E\) generated by the family of linear
functionals \(\mathcal F\). In particular, each \(\ell_{k,f,s}\in\mathcal F\)
is \(\tau\)-continuous (equivalently, the seminorms
\(p_{k,f,s}(\Phi)=|\ell_{k,f,s}(\Phi)|\) are \(\tau\)-continuous).

Conversely, for the weak topology \(\sigma(\mathcal E,\mathcal F)\), the
continuous dual is exactly \(\mathcal F\) (a standard fact for weak topologies
generated by a linear subspace of the algebraic dual). Hence every
\(\tau\)-continuous linear functional on \(\mathcal E\) belongs to
\(\mathcal F\).
\end{proof}

\begin{lemma}\label{lem:pairing-separates}
The family $\{\ell_{k,f,s}\}$ separates points of $\mathcal E$.
Equivalently, the pairing between matrix tests and maps is separating on
$\mathcal E$.
\end{lemma}

\begin{proof}
Let $\Phi,\Psi\in\mathcal E$ and assume
\[
\ell_{k,f,s}(\Phi)=\ell_{k,f,s}(\Psi)
\qquad\text{for all matrix tests }(k,f,s).
\]
In particular, for $k=1$ we have
\[
f(\Phi(s))=f(\Psi(s))
\qquad\text{for all }f\in S(\mathscr T),\ s\in\mathscr S.
\]
Since states separate points of the unital $C^*$-algebra $\mathscr T$, it
follows that $\Phi(s)=\Psi(s)$ for all $s\in\mathscr S$, hence $\Phi=\Psi$.
\end{proof}

\begin{corollary}\label{co_hau}
The topology $\tau=\sigma(\mathcal E,\mathcal F)$ is Hausdorff.
\end{corollary}

\begin{proof}
By Lemma~\ref{lem:pairing-separates}, the subspace $\mathcal F\subseteq\mathcal E^*$
separates points of $\mathcal E$.
\end{proof}

\begin{proposition}[Collapse of the matrix-level test topology]
\label{prop:tau-collapse}
For each \(r\in\mathbb N\), let \(\tau_r\) denote the locally convex topology on
\(\mathcal E\) generated by the seminorms
\[
p_{k,f,s}(\Phi)=\big|f(\Phi_k(s))\big|
\]
with \(1\le k\le r\), \(f\in S(M_k(\mathscr T))\), and \(s\in M_k(\mathscr S)\).
Then, for every \(r\ge 1\),
\[
\tau_r=\tau_1=\tau.
\]

Equivalently, every matrix-level test seminorm \(p_{k,f,s}\) is
\(\tau_1\)-continuous.
\end{proposition}

\begin{proof}
It is enough to show that every generating seminorm \(p_{k,f,s}\) of \(\tau\) is
\(\tau_1\)-continuous. Fix \(k\in\mathbb N\), a matrix state
\(f\in S(M_k(\mathscr T))\), and an element \(s\in M_k(\mathscr S)\). We divide the proof into three steps.

\medskip
\noindent\textbf{Step 1: Decompose the matrix-level test into matrix coefficients.}

Write \(s=[s_{ij}]_{i,j=1}^k\) with \(s_{ij}\in\mathscr S\). Let
\(\{E_{ij}\}_{i,j=1}^k\) be the standard matrix units in \(M_k\). Using the
identification \(M_k(\mathscr T)\cong M_k\otimes \mathscr T\), every
\(Y=[y_{ij}]\in M_k(\mathscr T)\) can be written as
\[
Y=\sum_{i,j=1}^k E_{ij}\otimes y_{ij}.
\]
Define linear functionals \(g_{ij}\in \mathscr T^*\) by
\[
g_{ij}(t):=f(E_{ij}\otimes t),\qquad t\in\mathscr T.
\]
Then for every \(Y=[y_{ij}]\in M_k(\mathscr T)\),
\[
f(Y)=\sum_{i,j=1}^k g_{ij}(y_{ij}).
\]
Applying this to \(Y=\Phi_k(s)=[\Phi(s_{ij})]\) gives
\begin{equation}\label{eq:matrix-coeff-decomp}
\ell_{k,f,s}(\Phi)
=f(\Phi_k(s))
=\sum_{i,j=1}^k g_{ij}\big(\Phi(s_{ij})\big)
\qquad(\Phi\in\mathcal E).
\end{equation}

\medskip
\noindent\textbf{Step 2: The complex span of states equals \(\mathscr T^*\).}

We use the standard fact that every bounded linear functional on a unital
\(C^*\)-algebra is a complex linear combination of states.  For completeness,
we record a short proof.

\begin{claim}\label{clm:states-span-dual}
For every \(g\in\mathscr T^*\), there exist states
\(\varphi_1,\varphi_2,\varphi_3,\varphi_4\in S(\mathscr T)\) and scalars
\(a_1,a_2,a_3,a_4\ge 0\) such that
\[
g
=
a_1\varphi_1-a_2\varphi_2
+i\bigl(a_3\varphi_3-a_4\varphi_4\bigr).
\]
\end{claim}

\begin{proof}[Proof of Claim~\ref{clm:states-span-dual}]
Define \(g^\sharp\in\mathscr T^*\) by
\[
g^\sharp(t):=\overline{g(t^*)}\qquad (t\in\mathscr T),
\]
and set
\[
u:=\frac{g+g^\sharp}{2},\qquad
v:=\frac{g-g^\sharp}{2i}.
\]
Then \(u\) and \(v\) are selfadjoint functionals, and \(g=u+iv\).

By the Jordan decomposition theorem for bounded selfadjoint functionals on a
\(C^*\)-algebra, there exist positive functionals \(u_+,u_-,v_+,v_-\in\mathscr T_+^*\)
such that
\[
u=u_+-u_-,\qquad v=v_+-v_-.
\]
Each nonzero positive functional \(h\in\mathscr T_+^*\) is a nonnegative scalar
multiple of a state: if \(a:=h(\1_{\mathscr T})>0\), then \(h=a\varphi\) with
\(\varphi:=a^{-1}h\in S(\mathscr T)\); if \(a=0\), then \(h=0\).

Applying this to \(u_\pm\) and \(v_\pm\) yields the desired representation of
\(g=u+iv\).
\end{proof}

\medskip
\noindent\textbf{Step 3: Dominate \(p_{k,f,s}\) by finitely many level-1 seminorms.}

Apply Claim~\ref{clm:states-span-dual} to each coefficient functional
\(g_{ij}\in\mathscr T^*\) from Step~1. Thus for each pair \((i,j)\) there exist
states \(\varphi_{ij}^{(1)},\dots,\varphi_{ij}^{(4)}\in S(\mathscr T)\) and
scalars \(a_{ij}^{(1)},\dots,a_{ij}^{(4)}\ge 0\) such that
\[
g_{ij}
=
a_{ij}^{(1)}\varphi_{ij}^{(1)}
-a_{ij}^{(2)}\varphi_{ij}^{(2)}
+i\Bigl(a_{ij}^{(3)}\varphi_{ij}^{(3)}-a_{ij}^{(4)}\varphi_{ij}^{(4)}\Bigr).
\]
Hence, for every \(\Phi\in\mathcal E\),
\[
\big|g_{ij}(\Phi(s_{ij}))\big|
\le
\sum_{m=1}^4 a_{ij}^{(m)}\,\big|\varphi_{ij}^{(m)}(\Phi(s_{ij}))\big|
=
\sum_{m=1}^4 a_{ij}^{(m)}\,p_{1,\varphi_{ij}^{(m)},s_{ij}}(\Phi).
\]
Combining this with \eqref{eq:matrix-coeff-decomp} and the triangle inequality,
we obtain
\begin{align*}
p_{k,f,s}(\Phi)
&=\big|f(\Phi_k(s))\big|
 =\left|\sum_{i,j=1}^k g_{ij}\big(\Phi(s_{ij})\big)\right| \\
&\le \sum_{i,j=1}^k \big|g_{ij}\big(\Phi(s_{ij})\big)\big| \\
&\le \sum_{i,j=1}^k\sum_{m=1}^4 a_{ij}^{(m)}\,
      p_{1,\varphi_{ij}^{(m)},s_{ij}}(\Phi).
\end{align*}
The right-hand side is a finite linear combination of seminorms generating
\(\tau_1\). Therefore \(p_{k,f,s}\) is \(\tau_1\)-continuous.

Since \(p_{k,f,s}\) was an arbitrary generating seminorm for \(\tau\), it
follows that \(\tau\subseteq\tau_1\). The reverse inclusion
\(\tau_1\subseteq\tau\) is immediate from the definitions, since the family
defining \(\tau\) already contains all level-1 tests. Hence \(\tau=\tau_1\).

Finally, for each \(r\ge 1\), we have \(\tau_1\subseteq\tau_r\subseteq\tau\),
and therefore
\[
\tau_r=\tau_1=\tau.
\]
\end{proof}

\begin{remark}\label{rem:tau-collapse-meaning}
Proposition~\ref{prop:tau-collapse} shows that higher matrix levels are not
needed to generate the weak topology \(\tau\) on \(\mathcal E\). Their role in
the present theory is instead geometric and order-theoretic: they enter through
the matrix-test inequalities, complete positivity, and the compatibility of the
duality with \(C^*\)-convex combinations.
\end{remark}


\section{Folding and Separation}
\label{sec:folding-separation}
This section develops the two ingredients used in the proof of the main
separation formulation theorem and its consequences.

The first is a finite-dimensional ``folding'' mechanism, which lets us combine
finitely many matrix tests into a single test at a larger matrix level.
The second is a separation lemma in the locally convex space
\((\mathcal E,\tau)\).  The separation step itself is standard
(Hahn--Banach), but we state it in a form adapted to our pairing between
completely positive maps and matrix tests.

\medskip

The following folding statement is the key input for the separation argument:
a finite complex linear combination of test functionals can be represented by a
single test functional (up to a positive scalar).  The phases are absorbed into
the test element, while the moduli become the coefficients of a convex
combination of matrix states.

\begin{lemma}\label{lem:folding-linear-combination}
Let $(k_j,f_j,s_j)$, $j=1,\dots,N$, be matrix tests, and let
$c_1,\dots,c_N\in\mathbb C$ be complex numbers, not all zero.
Write
\[
c_j=r_j e^{i\theta_j}\qquad (r_j\ge 0,\ \theta_j\in\mathbb R),
\]
and set $R:=\sum_{j=1}^N r_j>0$.
Define
\[
   K := \sum_{j=1}^N k_j,\qquad
   S := \diag\big(e^{i\theta_1}s_1,\dots,e^{i\theta_N}s_N\big)
      \in M_K(\mathscr{S}),
\]
and
\[
   F := \sum_{j=1}^N \lambda_j\, (f_j\circ\pi_j),
   \qquad \lambda_j := \frac{r_j}{R}\in[0,1],\quad \sum_{j=1}^N\lambda_j=1,
\]
where $\pi_j:M_K(\mathscr T)\to M_{k_j}(\mathscr T)$ extracts the $j$th
diagonal block.
Then $(K,F,S)$ is a matrix test and, for every $\Phi\in\mathcal E$,
\begin{equation}\label{eq:folding-linear-combination}
   \ell_{K,F,S}(\Phi)
   \;=\; \frac{1}{R}\sum_{j=1}^N c_j\, \ell_{k_j,f_j,s_j}(\Phi).
\end{equation}
In terms of the defining evaluations,
\eqref{eq:folding-linear-combination} reads
\[
   F\!\big(\Phi_K(S)\big)
   \;=\; \frac{1}{R}\sum_{j=1}^N c_j\,
   f_j\!\big(\Phi_{k_j}(s_j)\big)
   \qquad (\Phi\in\mathcal E).
\]
\end{lemma}

\begin{proof}
As above, each $f_j\circ\pi_j$ is a state on $M_K(\mathscr T)$, and hence so is
their convex combination $F$.

Fix $\Phi\in\mathcal E$. Since $\Phi$ is linear, its ampliation $\Phi_K$ is
linear on $M_K(\mathscr S)$, and block diagonality gives
\[
\Phi_K(S)=\diag\big(e^{i\theta_1}\Phi_{k_1}(s_1),\dots,e^{i\theta_N}\Phi_{k_N}(s_N)\big).
\]
Applying $F$ yields
\[
\ell_{K,F,S}(\Phi)
=F\!\big(\Phi_K(S)\big)
=\sum_{j=1}^N \lambda_j\, e^{i\theta_j}\, f_j\!\big(\Phi_{k_j}(s_j)\big)
=\frac1R\sum_{j=1}^N r_j e^{i\theta_j} f_j\!\big(\Phi_{k_j}(s_j)\big),
\]
which is exactly \eqref{eq:folding-linear-combination}.
\end{proof}

\medskip

We now turn to separation in the locally convex space $(\mathcal E,\tau)$,
where $\tau=\sigma(\mathcal E,\mathcal F)$ and
\[
\mathcal F=\mathrm{span}\{\ell_{k,f,s}\}
\]
(see Section~2). By Corollary~\ref{co_hau}, $(\mathcal E,\tau)$ is Hausdorff.

\begin{lemma}\label{lem:real_sep}
Let $\mathcal C\subseteq\mathcal E$ be $\tau$-closed and convex, and let
$\Phi_0\in\mathcal E\setminus\mathcal C$.
Then there exist a matrix test $(K,F,S)$ and a real number $\alpha$ such that
\begin{equation}\label{eq:real-separation}
   \Re\, \ell_{K,F,S}(\Phi_0) > \alpha
   \quad\text{and}\quad
   \Re\, \ell_{K,F,S}(\Phi) \le \alpha
   \quad \forall\,\Phi\in\mathcal C.
\end{equation}
In particular, if $\Phi_0\in\CP(\mathscr S,\mathscr T)$ and
$\mathcal C\subseteq\CP(\mathscr S,\mathscr T)$, then
\eqref{eq:real-separation} is equivalent to
\[
   \Re\, F\!\big(\Phi_{0,K}(S)\big) > \alpha
   \quad\text{and}\quad
   \Re\, F\!\big(\Phi_K(S)\big) \le \alpha
   \quad \forall\,\Phi\in\mathcal C.
\]
\end{lemma}

\begin{proof}
By the Hahn--Banach separation theorem in the Hausdorff locally convex space
$(\mathcal E,\tau)$ (viewed as a real locally convex space), there exist a
$\tau$-continuous real-linear functional $L:\mathcal E\to\mathbb R$ and
$\alpha\in\mathbb R$ such that
\[
L(\Phi_0)>\alpha\ge L(\Phi)\qquad(\Phi\in\mathcal C).
\]

Define a complex-linear functional $\widetilde L:\mathcal E\to\mathbb C$ by
\[
\widetilde L(\Phi):=L(\Phi)-iL(i\Phi).
\]
Then $\widetilde L$ is $\tau$-continuous and
\[
L(\Phi)=\Re\,\widetilde L(\Phi)\qquad(\Phi\in\mathcal E).
\]
By Lemma~\ref{lem:tau-dual}, $\widetilde L\in\mathcal F$. Since
$\mathcal F=\mathrm{span}\{\ell_{k,f,s}\}$ (algebraic span), there exist
$N\in\mathbb N$, matrix tests $(k_j,f_j,s_j)$, and scalars $c_j\in\mathbb C$
such that
\[
\widetilde L=\sum_{j=1}^N c_j\,\ell_{k_j,f_j,s_j}.
\]
Discarding zero coefficients if necessary, we may assume the $c_j$ are not all
zero. Applying Lemma~\ref{lem:folding-linear-combination}, we obtain a single
matrix test $(K,F,S)$ and a scalar $R>0$ such that
\[
\widetilde L(\Phi)=R\,\ell_{K,F,S}(\Phi)\qquad(\Phi\in\mathcal E).
\]
Therefore
\[
L(\Phi)=\Re\,\widetilde L(\Phi)=R\,\Re\,\ell_{K,F,S}(\Phi)
\qquad(\Phi\in\mathcal E).
\]
Substituting this into the separating inequalities for $L$ gives
\[
\Re\,\ell_{K,F,S}(\Phi_0)>\alpha/R
\quad\text{and}\quad
\Re\,\ell_{K,F,S}(\Phi)\le\alpha/R\quad(\Phi\in\mathcal C).
\]
Renaming $\alpha/R$ as $\alpha$ proves \eqref{eq:real-separation}.
\end{proof}

\section{Support-Function Characterization and Consequences}
\label{sec:matrix-bipolar}
 This section contains the main geometric consequence of the matrix-test
framework. We prove a support-function/separation characterization of the
$\tau$-closed $C^*$-convex hull of a family of completely positive maps, and
then derive several concrete consequences.

We continue to use the one-sided polar convention. Accordingly, polar
inequalities are written in terms of real parts:
\[
\Re\,\langle (k,f,s),\Phi\rangle \le 1.
\]
This convention is well suited to convex separation and leads to a clean
support-function description of $\overline{\cconv(\mathcal K)}^{\,\tau}$.

We begin by showing that complete positivity is closed for the weak topology
$\tau=\sigma(\mathcal E,\mathcal F)$ introduced in Section~2.

\begin{lemma}\label{lem:CP-tau-closed}
The set $\CP(\mathscr S,\mathscr T)$ is $\tau$-closed in $\mathcal E$.
\end{lemma}

\begin{proof}
Let $(\Phi_\nu)$ be a net in $\CP(\mathscr S,\mathscr T)$ converging to
$\Phi\in\mathcal E$ in the topology $\tau$.
We show that $\Phi$ is completely positive.

Fix $k\in\mathbb N$ and $x\in M_k(\mathscr S)_+$. For every state
$f\in S(M_k(\mathscr T))$, each $\Phi_\nu$ is completely positive, so
\[
f\!\big((\Phi_\nu)_k(x)\big)\ge 0.
\]
Since $\Phi_\nu\to\Phi$ in $\tau$, we have
\[
\ell_{k,f,x}(\Phi_\nu)\to \ell_{k,f,x}(\Phi),
\]
that is,
\[
f\!\big((\Phi_\nu)_k(x)\big)\to f\!\big(\Phi_k(x)\big).
\]
Hence $f(\Phi_k(x))\ge 0$ for every state $f$ on $M_k(\mathscr T)$.

By the standard positivity criterion in a unital $C^*$-algebra
(an element $a$ is positive iff $f(a)\ge 0$ for all states $f$),
it follows that $\Phi_k(x)\ge 0$ in $M_k(\mathscr T)$.
Since $k$ and $x\in M_k(\mathscr S)_+$ were arbitrary, $\Phi$ is completely
positive.
\end{proof}

\medskip

For $\mathcal K\subseteq\CP(\mathscr S,\mathscr T)$ and a matrix test
$(k,f,s)$, define the support function
\begin{equation}\label{eq:support-function}
h_{\mathcal K}(k,f,s)
:=
\sup_{\Psi\in\cconv(\mathcal K)} \Re\, f\!\big(\Psi_k(s)\big)
\in [-\infty,\infty].
\end{equation}
(If $\mathcal K=\varnothing$, then $\cconv(\mathcal K)=\varnothing$ and the
supremum is $-\infty$.)

The next theorem is a convex-duality statement in the weak topology
$\tau=\sigma(\mathcal E,\mathcal F)$. What is specific to the present setting
is that, after the folding results of Section~4, the relevant
$\tau$-continuous functionals can be represented by matrix tests. This turns
the abstract separation criterion into a concrete family of matrix-test
inequalities.

\begin{theorem}[Separation formulation]\label{thm:separation-formulation}
Let $\mathscr S$ be an operator system, let $\mathscr T$ be a unital
$C^*$-algebra, let $\mathcal K\subseteq\CP(\mathscr S,\mathscr T)$, and let
$\Phi_0\in\CP(\mathscr S,\mathscr T)$.
Then the following are equivalent:
\begin{enumerate}
\item[(i)] $\displaystyle \Phi_0\in \overline{\cconv(\mathcal K)}^{\,\tau}$.
\item[(ii)] For every matrix test $(k,f,s)$,
\begin{equation}\label{eq:separation-formulation}
\Re\, f\!\big((\Phi_0)_k(s)\big)
\le
h_{\mathcal K}(k,f,s)
=
\sup_{\Psi\in\cconv(\mathcal K)}
\Re\, f\!\big(\Psi_k(s)\big).
\end{equation}
\end{enumerate}
\end{theorem}

\begin{proof}
Set
\[
C:=\overline{\cconv(\mathcal K)}^{\,\tau}\subseteq \mathcal E.
\]
By Lemma~\ref{lem:CP-tau-closed}, we have
$C\subseteq \CP(\mathscr S,\mathscr T)$.
Moreover, $\cconv(\mathcal K)$ is (ordinarily) convex (scalar convex
combinations are special $C^*$-convex combinations), hence its $\tau$-closure
$C$ is convex.

Assume first that $\Phi_0\in C$. Let $(k,f,s)$ be any matrix test. Since
$\ell_{k,f,s}$ is $\tau$-continuous (Lemma~\ref{lem:tau-dual}), the map
$\Phi\mapsto \Re\,\ell_{k,f,s}(\Phi)$ is $\tau$-continuous on $\mathcal E$.
Therefore
\[
\Re\,\ell_{k,f,s}(\Phi_0)
\le
\sup_{\Psi\in C}\Re\,\ell_{k,f,s}(\Psi).
\]
Because $\cconv(\mathcal K)$ is $\tau$-dense in $C$, and
$\Re\,\ell_{k,f,s}$ is continuous, the supremum over $C$ agrees with the
supremum over $\cconv(\mathcal K)$ (possibly equal to $+\infty$). Hence
\[
\Re\, f\!\big((\Phi_0)_k(s)\big)
=
\Re\,\ell_{k,f,s}(\Phi_0)
\le
\sup_{\Psi\in\cconv(\mathcal K)} \Re\,\ell_{k,f,s}(\Psi),
\]
which is exactly \eqref{eq:separation-formulation}.

Conversely, assume (ii) holds and suppose, for contradiction, that
$\Phi_0\notin C$. Since $C$ is $\tau$-closed and convex,
Lemma~\ref{lem:real_sep} applies and yields a matrix test $(K,F,S)$ and
$\alpha\in\mathbb R$ such that
\[
\Re\,\ell_{K,F,S}(\Phi_0)>\alpha
\quad\text{and}\quad
\Re\,\ell_{K,F,S}(\Phi)\le\alpha \qquad(\Phi\in C).
\]
In particular,
\[
\sup_{\Psi\in\cconv(\mathcal K)} \Re\,\ell_{K,F,S}(\Psi)\le \alpha
< \Re\,\ell_{K,F,S}(\Phi_0),
\]
which contradicts (ii) for the test $(K,F,S)$.
Therefore $\Phi_0\in C$, proving (i).
\end{proof}

\begin{corollary}[Single-test witness of non-membership]
\label{cor:single-test-witness}
Let $\mathcal K\subseteq\CP(\mathscr S,\mathscr T)$ and let
$\Phi_0\in\CP(\mathscr S,\mathscr T)$.
If
\[
\Phi_0\notin \overline{\cconv(\mathcal K)}^{\,\tau},
\]
then there exists a matrix test $(K,F,S)$ such that
\[
\Re\,F\!\big((\Phi_0)_K(S)\big)
>
\sup_{\Psi\in\cconv(\mathcal K)} \Re\,F\!\big(\Psi_K(S)\big).
\]
\end{corollary}

\begin{proof}
If $\Phi_0\notin \overline{\cconv(\mathcal K)}^{\,\tau}$, then condition
\textup{(i)} in Theorem~\ref{thm:separation-formulation} fails. Hence
condition \textup{(ii)} fails as well. Therefore there exists a matrix test
$(K,F,S)$ such that
\[
\Re\,F\!\big((\Phi_0)_K(S)\big)
>
\sup_{\Psi\in\cconv(\mathcal K)} \Re\,F\!\big(\Psi_K(S)\big).
\]
\end{proof}

\begin{corollary}[Support-function criterion for inclusion]
\label{cor:support-inclusion}
Let $\mathcal K,\mathcal L\subseteq\CP(\mathscr S,\mathscr T)$. Then the
following are equivalent:
\begin{enumerate}
\item[(i)]
\[
\overline{\cconv(\mathcal K)}^{\,\tau}
\subseteq
\overline{\cconv(\mathcal L)}^{\,\tau}.
\]
\item[(ii)] For every matrix test $(k,f,s)$,
\[
h_{\mathcal K}(k,f,s)\le h_{\mathcal L}(k,f,s).
\]
\end{enumerate}
\end{corollary}

\begin{proof}
Assume (i). Let $(k,f,s)$ be a matrix test. Since
\[
\cconv(\mathcal K)\subseteq \overline{\cconv(\mathcal K)}^{\,\tau}
\subseteq \overline{\cconv(\mathcal L)}^{\,\tau},
\]
we have
\[
h_{\mathcal K}(k,f,s)
=
\sup_{\Psi\in\cconv(\mathcal K)} \Re\,f\!\big(\Psi_k(s)\big)
\le
\sup_{\Phi\in \overline{\cconv(\mathcal L)}^{\,\tau}}
\Re\,f\!\big(\Phi_k(s)\big).
\]
Because $\Re\,\ell_{k,f,s}$ is $\tau$-continuous and
$\cconv(\mathcal L)$ is $\tau$-dense in its closure, the supremum over
$\overline{\cconv(\mathcal L)}^{\,\tau}$ agrees with the supremum over
$\cconv(\mathcal L)$. Hence
\[
h_{\mathcal K}(k,f,s)\le h_{\mathcal L}(k,f,s).
\]

Conversely, assume (ii), and let
\[
\Phi_0\in \overline{\cconv(\mathcal K)}^{\,\tau}.
\]
By Lemma~\ref{lem:CP-tau-closed}, we have \(\Phi_0\in\CP(\mathscr S,\mathscr T)\).
By Theorem~\ref{thm:separation-formulation}, for every matrix test $(k,f,s)$,
\[
\Re\,f\!\big((\Phi_0)_k(s)\big)\le h_{\mathcal K}(k,f,s)\le h_{\mathcal L}(k,f,s).
\]
Applying Theorem~\ref{thm:separation-formulation} again (with \(\mathcal L\) in
place of \(\mathcal K\)), we conclude that
\[
\Phi_0\in \overline{\cconv(\mathcal L)}^{\,\tau}.
\]
Thus
\[
\overline{\cconv(\mathcal K)}^{\,\tau}
\subseteq
\overline{\cconv(\mathcal L)}^{\,\tau}.
\]
\end{proof}

\begin{corollary}[Support-function reconstruction]
\label{cor:support-reconstruction}
Let $\mathcal K,\mathcal L\subseteq\CP(\mathscr S,\mathscr T)$. Then
\[
\overline{\cconv(\mathcal K)}^{\,\tau}
=
\overline{\cconv(\mathcal L)}^{\,\tau}
\]
if and only if
\[
h_{\mathcal K}(k,f,s)=h_{\mathcal L}(k,f,s)
\qquad\text{for every matrix test }(k,f,s).
\]
\end{corollary}

\begin{proof}
This follows immediately from Corollary~\ref{cor:support-inclusion} by applying
it in both directions.
\end{proof}

\medskip

We now return to the normalized one-sided polar. In the notation of
Section~2 (with the one-sided polar convention),
\[
\mathcal K^\circ
=
\Bigl\{(k,f,s)\ \Big|\
\Re\, f\!\big(\Psi_k(s)\big)\le 1
\ \text{for all }\Psi\in\cconv(\mathcal K)\Bigr\},
\]
and for a collection of matrix tests $\mathcal C$,
\[
\mathcal C^{\circ_{C^*}}
=
\Bigl\{\Phi\in\CP(\mathscr S,\mathscr T)\ \Big|\
\Re\, f\!\big(\Phi_k(s)\big)\le 1
\ \text{for all }(k,f,s)\in\mathcal C\Bigr\}.
\]

The next corollary shows that under the natural additional hypothesis
\(0\in \overline{\cconv(\mathcal K)}^{\,\tau}\), the normalized bipolar does
recover the closed $C^*$-convex hull exactly.

\begin{corollary}\label{cor:exact-bipolar-if-zero}
If \(0\in \overline{\cconv(\mathcal K)}^{\,\tau}\), then
\[
(\mathcal K^\circ)^{\circ_{C^*}}
=
\overline{\cconv(\mathcal K)}^{\,\tau}.
\]
\end{corollary}

\begin{proof}
 
Set
\[
C:=\overline{\cconv(\mathcal K)}^{\,\tau}
\subseteq \CP(\mathscr S,\mathscr T).
\]
We prove
\[
C=(\mathcal K^\circ)^{\circ_{C^*}}
\]
under the assumption \(0\in C\).

First, \(C\subseteq(\mathcal K^\circ)^{\circ_{C^*}}\).
Indeed, let \(\Phi\in C\) and let
\((k,f,s)\in\mathcal K^\circ\). By the definition of
\(\mathcal K^\circ\),
\[
\sup_{\Psi\in\cconv(\mathcal K)}
\Re\,f\bigl(\Psi_k(s)\bigr)\le 1.
\]
Since \(\Phi\in \overline{\cconv(\mathcal K)}^{\,\tau}\) and
\(\Phi\mapsto \Re\,f(\Phi_k(s))\) is \(\tau\)-continuous, it follows that
\[
\Re\, f(\Phi_k(s))\le 1.
\]
Thus \(\Phi\in (\mathcal K^{\circ_{C^*}})^{\circ_{C^*}}\).

For the reverse inclusion, let
\[
\Phi_0\in(\mathcal K^\circ)^{\circ_{C^*}}.
\]
We will verify the inequalities in
Theorem~\ref{thm:separation-formulation}, which will imply \(\Phi_0\in C\).

Fix a matrix test \((k,f,s)\), and write
\[
h:=h_{\mathcal K}(k,f,s)
=\sup_{\Psi\in\cconv(\mathcal K)} \Re\, f(\Psi_k(s))\in [ -\infty,\infty ].
\]
Because \(0\in C\) and \(\Phi\mapsto \Re\,f(\Phi_k(s))\) is \(\tau\)-continuous,
we have
\[
h=h_C(k,f,s)\ge \Re\,f(0)=0.
\]
Hence \(h\in[0,\infty]\).

We now consider three cases.

\smallskip
\noindent\emph{Case 1: \(h=\infty\).}
Then the desired inequality
\[
\Re\,f((\Phi_0)_k(s))\le h_{\mathcal K}(k,f,s)
\]
is automatic.

\smallskip
\noindent\emph{Case 2: \(0<h<\infty\).}
By positive homogeneity of the support function in the test variable \(s\),
\[
h_{\mathcal K}\!\left(k,f,\frac{s}{h}\right)=\frac{1}{h}\,h_{\mathcal K}(k,f,s)=1.
\]
Therefore $\left(k,f,\frac{s}{h}\right)\in\mathcal K^\circ$.  Since
\(\Phi_0\in(\mathcal K^\circ)^{\circ_{C^*}}\), we get
\[
\Re\, f\!\left((\Phi_0)_k\!\left(\frac{s}{h}\right)\right)\le 1.
\]
Multiplying by \(h\) yields
\[
\Re\, f((\Phi_0)_k(s))\le h=h_{\mathcal K}(k,f,s).
\]

\smallskip
\noindent\emph{Case 3: \(h=0\).}
For every \(t>0\), again by positive homogeneity,
\[
h_{\mathcal K}(k,f,ts)=t\,h_{\mathcal K}(k,f,s)=0\le 1,
\]
so \((k,f,ts)\in\mathcal K^\circ\).  Since
\(\Phi_0\in(\mathcal K^\circ)^{\circ_{C^*}}\), we have
\[
\Re\, f((\Phi_0)_k(ts))\le 1.
\]
By linearity of \((\Phi_0)_k\),
\[
t\,\Re\, f((\Phi_0)_k(s))\le 1
\quad\text{for all }t>0.
\]
Dividing by \(t\) and letting \(t\to\infty\), we obtain
\[
\Re\, f((\Phi_0)_k(s))\le 0=h_{\mathcal K}(k,f,s).
\]

Thus in all cases,
\[
\Re\, f((\Phi_0)_k(s))\le h_{\mathcal K}(k,f,s)
\]
for every matrix test \((k,f,s)\).  By
Theorem~\ref{thm:separation-formulation}, it follows that \(\Phi_0\in C\).

Therefore
\[
C=(\mathcal K^\circ)^{\circ_{C^*}},
\]
as claimed.
\end{proof}


\section*{Acknowledgment}
The authors thank the anonymous referee for a careful reading of the
manuscript and for constructive comments that helped improve its
presentation and clarify its relation to matrix and noncommutative
convexity.

The research of the second author was funded by the European Union’s NextGenerationEU programme, as part of the institutional project “Resilient Self-Healing Future Power Systems -- RePowerFER,” which is included in the programme agreement of the University of Zagreb Faculty of Electrical Engineering and Computing. The views and opinions expressed are solely those of the author and do not necessarily reflect the official position of the European Union or the European Commission. Neither the European Union nor the European Commission can be held responsible for them.

\end{document}